\documentclass[twoside,a4paper]{article} 

\usepackage[margin=2cm]{geometry}
\usepackage{amsmath,amsthm,amssymb}

\newcommand{\wap}{\operatorname{WAP}}
\newcommand{\proten}{\widehat\otimes}
\newcommand{\mc}{\mathcal}
\newcommand{\ip}[2]{\langle #1 , #2 \rangle}

\theoremstyle{plain}
\newtheorem{proposition}{Proposition}[section]
\newtheorem{theorem}[proposition]{Theorem}
\newtheorem{corollary}[proposition]{Corollary}
\newtheorem{lemma}[proposition]{Lemma}
\theoremstyle{definition}

\theoremstyle{remark}

\begin{document}

\large
\title{\textsc{A bicommutant theorem for dual Banach algebras}}
\author{Matthew Daws}
\maketitle

\begin{abstract}
A dual Banach algebra is a Banach algebra which is a dual space, with the
multiplication being separately weak$^*$-continuous.  We show that given a unital
dual Banach algebra $\mc A$, we can find a reflexive Banach space $E$, and an
isometric, weak$^*$-weak$^*$-continuous homomorphism $\pi:\mc A\rightarrow\mc B(E)$
such that $\pi(\mc A)$ equals its own bicommutant.


2000/2010 \emph{Mathematical Subject Classification:}
46H05, 46H15, 47L10 (primary), 46A32, 46B10
\end{abstract}

\section{Introduction}

Given a Banach space $E$, we write $\mc B(E)$ for the Banach algebra of operators
on $E$.  Given a subset $X\subseteq\mc B(E)$, we write $X'$ for the commutant of $X$,
\[ X' = \{ T\in\mc B(E) : TS=ST \ (S\in X) \}. \]
The von Neumann bicommutant theorem tells us that if $E$ is a Hilbert space,
and $X$ is a $*$-closed, unital subalgebra, then $X''$ is the strong operator
topology closure of $X$ in $\mc B(E)$.  If $X$ is not $*$-closed, then this result
may fail (consider strictly upper-triangular two-by-two matricies).  However, a result
of Blecher and Solel, \cite{BS}, shows, in particular, that if $X$ is weak$^*$-closed,
that we can find another Hilbert space $K$, and a completely isometric,
weak$^*$-weak$^*$-continuous homomorphism $\pi:X\rightarrow \mc B(K)$, such that
$\pi(X) = \pi(X)''$.  That is, if we change the Hilbert space which our algebra
acts on, we do have a bicommutant theorem.

A dual Banach algebra is a Banach algebra which is a dual space, such that the
multiplication is weak$^*$-continuous.  Building on work of Young and Kaiser, the
author showed in \cite{Daws} that given a dual Banach algebra $\mc A$, we can find
a reflexive Banach space $E$ and an isometric, weak$^*$-weak$^*$-continuous
homomorphism $\pi:\mc A\rightarrow\mc B(E)$.  In this paper, we show that when
$\mc A$ is unital, we can choose $E$ and $\pi$ such that $\pi(\mc A) = \pi(\mc A)''$.
The method is similar to that used in \cite{BS} (although we follow the presentation
of \cite{BLM}) combined with an idea adapted from \cite[Section~6]{Daws}.

\subsection{Acknowledgments}

The author wishes to thank Stuart White and Allan Sinclair for suggesting this
problem, and to thank David Blecher for bringing \cite{BS} to his attention.

\section{Notation and preliminary results}

Given a Banach space $E$, let $E^*$ be the dual space to $E$.  For $\mu\in E^*$
and $x\in E$, we write $\ip{\mu}{x} = \mu(x)$.  For $X\subseteq E$, let
\[ X^\perp = \{ \mu\in E^* : \ip{\mu}{x}=0 \ (x\in X) \}. \]
For $Y\subseteq E^*$, let
\[ {}^\perp Y = \{ x\in E : \ip{\mu}{x}=0  \ (\mu\in Y) \}. \]
Then ${}^\perp(X^\perp)$ is the closure of the linear span of $X$, while
$({}^\perp Y)^\perp$ is the weak$^*$-closure of the linear span of $Y$.
We may canonically identify $X^*$ with $E/X^\perp$, and $(E/X)^*$ with
$X^\perp$.  In particular, $Y$ is weak$^*$-closed if and only if
$Y = ({}^\perp Y)^\perp$, and in this case, the canonical predual of $Y$
is $E / {}^\perp Y$.

We write $E^*\proten E$ for the projective tensor product of $E^*$ with $E$.  This is
the completion of the algebraic tensor product $E^*\otimes E$ with respect to the norm
\[ \|\tau\|_\pi = \inf\Big\{ \sum_{k=1}^n \|\mu_k\| \|x_k\| :
\tau = \sum_{k=1}^n \mu_k\otimes x_k \Big\}. \]
Any element of $E^*\proten E$ can be written as $\sum_k \mu_k\otimes x_k$ with
$\sum_k \|\mu_k\| \|x_k\|<\infty$.  For further details, see \cite{Dales} or
\cite{PalBook}, for example.

The Banach algebra $\mc B(E)$ is a dual Banach algebra with respect to the
predual $E^*\proten E$, the dual pairing being given by
\[ \ip{T}{\mu\otimes x} = \ip{\mu}{T(x)}
\qquad (T\in\mc B(E), \mu\otimes x\in E^*\proten E), \]
and linearity and continuity.  Indeed, under many circumstances, this is the
unique predual for $\mc B(E)$, see \cite[Theorem~4.4]{Daws}.

It follows that any weak$^*$-closed subalgebra of $\mc B(E)$ is also a dual
Banach algebra: then \cite[Corollary~3.8]{Daws} shows that every dual Banach
algebra arises in this way.  If $X\subseteq\mc B(E)$, then $X'$ is a closed
subalgebra of $\mc B(E)$.  Notice that $T\in X'$ if and only if $T$
annihilates all $\tau\in E^*\proten E$ of the form
\[ \tau = \mu\otimes S(x) - S^*(\mu)\otimes x
\qquad (S\in X, \mu\in E^*, x\in E). \]
Hence $X' = Y^\perp = (E^*\proten E / Y)^*$ is weak$^*$-closed, where
$Y$ is the closed linear span of such $\tau$.  In particular, $X''$ is
a weak$^*$-closed subalgebra of $\mc B(E)$ containing $X$, and so $X''$
contains the weak$^*$-closed algebra generated by $X$.

We shall follow the ideas of \cite[Theorem~3.2.14]{BLM}; see \cite{BS}
for a fuller treatment.  We first establish some preliminary results.
Given a Banach space $E$, we write $\ell^2(E)$ for the Banach space
consisting of sequences $(x_n)$ in $E$ with norm $\|(x_n)\|_2 = 
\Big( \sum_n \|x_n\|^2 \Big)^{1/2}$.  Throughout, we could instead work
with $\ell^p(E)$ for $1<p<\infty$, if we so wished.
Then $\ell^2(E)^* = \ell^2(E^*)$,
and $\ell^2(E)$ is reflexive if $E$ is.  For each $n$, let $\iota_n:
E\rightarrow\ell^2(E)$ be the injection onto the $n$th co-ordinate, and
let $P_n:\ell^2(E)\rightarrow E$ be the projection onto the $n$th co-ordinate.
For $T\in\mc B(E)$, let $T^{(\infty)}\in\mc B(\ell^2(E))$ be the operator
given by applying $T$ to each co-ordinate.  Notice that $T^{(\infty)} \iota_n
= \iota_n T$ and $P_n T^{(\infty)} = T P_n$, for each $n$.
For $X\subseteq\mc B(E)$,
let $X^{(\infty)} = \{ T^{(\infty)} : T\in X \}$.  Given a homomorphism
$\pi:\mc A\rightarrow\mc B(E)$, let $\pi^{(\infty)}:\mc A
\rightarrow\mc B(\ell^2(E))$ by the homomorphism given by $\pi^{(\infty)}(a)
= \pi(a)^{(\infty)}$ for each $a\in\mc A$.

\begin{lemma}\label{lemma::one}
For a Banach space $E$, and $X\subseteq\mc B(E)$, we have that
$(X^{(\infty)})'' = (X'')^{(\infty)}$.
\end{lemma}
\begin{proof}
Let $Q\in(X^{(\infty)})'$.  For $n,m\in\mathbb N$ and $S\in X$,
we have that $P_n Q \iota_m S = P_n Q S^{(\infty)} \iota_m
= P_n S^{(\infty)} Q \iota_m = S P_n Q \iota_m$.
Thus $P_n Q^{(\infty)} \iota_m \in X'$, for each $n,m$.  Similarly,
one can show that for $Q\in\mc B(\ell^2(E))$, if $P_n Q \iota_m
\in X'$ for all $n,m$, then $Q\in(X^{(\infty)})'$.

So, given $T\in X''$ and $Q\in(X^{(\infty)})'$, we have that $TP_nQ\iota_m
= P_nQ\iota_mT$ for all $n,m$.  Thus, for all $n,m$, it follows that
$P_n T^{(\infty)} Q \iota_m = P_n Q T^{(\infty)} \iota_m$, from which it
follows that $T^{(\infty)} Q = Q T^{(\infty)}$.  Thus $(X'')^{(\infty)}
\subseteq (X^{(\infty)})''$.

For the converse, let $T \in (X^{(\infty)})''$.  For each $n,m$, notice that
$\iota_n P_m \in (X^{(\infty)})'$, so that $T \iota_n P_m = \iota_n P_m T$.
Let $r\in\mathbb N$, so that
\[ T \iota_n \delta_{m,r} = T \iota_n P_m \iota_r = \iota_n P_m T \iota_r. \]
It follows that $T\iota_r = \iota_r R$ for some $R\in\mc B(E)$, and that $R$ does
not depend upon $r$.  Thus there must exist $R\in\mc B(E)$ with $T=R^{(\infty)}$.
Now let $S\in X'$, so that $S^{(\infty)} \in (X^{(\infty)})'$, and hence
\[ (RS)^{(\infty)} = T S^{(\infty)} = S^{(\infty)} T = (SR)^{(\infty)}. \]
It follows that $R\in X''$, and hence that $(X^{(\infty)})'' \subseteq
(X'')^{(\infty)}$.
\end{proof}

\begin{lemma}\label{lemma::two}
Let $E$ be a reflexive Banach space, and let $X\subseteq\mc B(E)$ be a subalgebra.
Let $X_w$ be the weak$^*$-closure of $X$ in $\mc B(E)$, with respect to the
predual $E^*\proten E$.  Then $(X_w)^{(\infty)} = (X^{(\infty)})_w$.
\end{lemma}
\begin{proof}
Let $T\in (X^{(\infty)})_w$.  For $x\in E,\mu\in E^*$ and $n\not=m$, certainly
$\iota_n(\mu) \otimes \iota_m(x) \in {}^\perp(X^{(\infty)})$, and so
\[ 0 = \ip{\iota_n(\mu)}{T\iota_m(x)} = \ip{\mu}{P_n T \iota_m(x)}. \]
Thus $P_n T \iota_m=0$ whenever $n\not=m$.  For any $x,\mu,n$ and $m$,
we also have that
\[ \iota_n(\mu)\otimes\iota_n(x) - \iota_m(\mu)\otimes\iota_m(x) \in
{}^\perp(X^{(\infty)}). \]
It follows that $P_n T \iota_n = P_m T \iota_m$.  Combining these results,
we conclude that $T=S^{(\infty)}$ for some $S\in\mc B(E)$.

Let $\tau \in {}^\perp X \subseteq E^*\proten E$, say $\tau = \sum_k \mu_k\otimes x_k$.
For $R\in X$ and each $n$, we have that
\[ \ip{R^{(\infty)}}{\sum_k \iota_n(\mu_k) \otimes \iota_n(x_k)} = 0, \]
so that $\sigma=\sum_k \iota_n(\mu_k) \otimes \iota_n(x_k) \in {}^\perp (X^{(\infty)})$.
So
\[ 0 = \ip{T}{\sigma} = \ip{S^{(\infty)}}{\sigma} = \ip{S}{\tau}, \]
from which it follows that $S\in X_w$.  So $(X^{(\infty)})_w \subseteq (X_w)^{(\infty)}$.

For the converse, let $T\in X_w$, and let $\tau\in {}^\perp (X^{(\infty)})$, say
$\tau = \sum_n \mu_n \otimes x_n$.  By rescaling, we may suppose that $\sum_n \|\mu_n\|^2
= \sum_n \|x_n\|^2 <\infty$.  For each $n$, we have that $\mu_n=(\mu^{(n)}_k)$, say,
where $\|\mu_n\|^2 = \sum_k \|\mu^{(n)}_k\|^2$.  Thus $\sum_{n,k} \|\mu^{(n)}_k\|^2<\infty$.
Similarly, each $x_n = (x^{(n)}_k)$, and $\sum_{n,k} \|x^{(n)}_k\|^2<\infty$.  We can
now compute that, for $S\in X$,
\[ 0 = \ip{S^{(\infty)}}{\tau} = \sum_n \ip{\mu_n}{S^{(\infty)}(x_n)}
= \sum_{n,k} \ip{\mu^{(n)}_k}{S(x^{(n)}_k)}, \]
so that $\sigma = \sum_{n,k} \mu^{(n)}_k \otimes x^{(n)}_k \in {}^\perp X$
(where this sum converges absolutely by an application of the Cauchy-Schwarz inequality).
Then $0 = \ip{T}{\sigma} = \ip{T^{(\infty)}}{\tau}$, from which it follows that
$T^{(\infty)} \in (X^{(\infty)})_w$.  So $(X_w)^{(\infty)} \subseteq (X^{(\infty)})_w$.
\end{proof}

The following lemma is usually stated in terms of ``reflexivity'' of a subspace
of $\mc B(E)$, but this is a different meaning to that of a reflexive Banach space,
so we avoid this terminology.

\begin{lemma}\label{lemma::three}
Let $E$ be a reflexive Banach space, and let $X\subseteq\mc B(E)$ be a weak$^*$-closed
subspace.  If $T\in\mc B(\ell^2(E))$ is such that, for each $x\in\ell^2(E)$,
we have that $T(x)$ is in the closure of $\{ S^{(\infty)}(x) : S\in X \}$, then
actually $T\in X^{(\infty)}$.
\end{lemma}
\begin{proof}
Let $T$ be as stated, so for each $n$, we have that the image of $T\iota_n$ is a
subset of the image of $\iota_n$.  By considering what $T$ maps $(\iota_1+\cdots+\iota_n)(x)$
to, for any $x\in E$, we may conclude that $T=R^{(\infty)}$ for some $R\in\mc B(E)$.

Let $\tau \in {}^\perp X$, say $\tau = \sum_n \mu_n\otimes x_n$, where we may suppose
that $\sum_n \|\mu_n\|^2 = \sum_n \|x_n\|^2 < \infty$.  Let $\mu=(\mu_n)\in\ell^2(E^*)$
and $x=(x_n)\in\ell^2(E)$, so that
\[ \ip{R}{\tau} = \ip{\mu}{R^{(\infty)}(x)} = \ip{\mu}{T(x)}. \]
However, notice that $\ip{\mu}{S^{(\infty)}(x)} = \ip{S}{\tau}=0$ for each $S\in X$,
so by the assumption on $T$, it follows also that $\ip{\mu}{T(x)}=0$, so $\ip{R}{\tau}=0$.
So $R\in ({}^\perp X)^\perp = X$, that is, $T\in X^{(\infty)}$.
\end{proof}

\section{The main result}

Let us introduce some temporary terminology, motivated by \cite{BLM}.  Let
$\mc A$ be a Banach algebra, and $E$ be a left $\mc A$-module (which we assume to
be a Banach space with contractive actions).  In this section,
we shall always suppose that $E$ is \emph{essential}, that is, the linear span of
$\{ a\cdot x: a\in\mc A, x\in E \}$ is dense in $E$.

We say that $E$ is \emph{cyclic} if there exists $x\in E$ with $\mc A\cdot x =
\{ a\cdot x: a\in\mc A\}$ being dense in $E$.  We say that $E$ is
\emph{self-generating} if, for each closed cyclic submodule $K\subseteq E$, the
linear span of $\{ T(E) : T:E\rightarrow K \text{ is an $\mc A$-module homomomorphism}\}$
is dense in $K$.

The following is very similar to the presentation in \cite{BLM}, but we check that the
details still work for reflexive Banach spaces, and not just Hilbert spaces.

\begin{theorem}\label{thm::one}
Let $\mc A$ be a unital Banach algebra, and let
$E$ be a reflexive Banach space with a bounded homomorphism $\pi:\mc A\rightarrow
\mc B(E)$.  Use $\pi$ to turn $E$ into a left $\mc A$-module, and suppose that
$\ell^2(E)$ is self-generating.  Then $\pi(\mc A)''$ agrees with the weak$^*$-closure
of $\pi(\mc A)$ in $\mc B(E)$.
\end{theorem}
\begin{proof}
Let $\mc B$ be the closure of $\pi(\mc A)$ in $\mc B(E)$, and let $\mc B_w$ be
the weak$^*$-closure of $\mc B$.  We wish to show that $\mc B_w = \mc B''$.

Let $T\in (\mc B'')^{(\infty)} \subseteq \mc B(\ell^2(E))$, let
$x\in\ell^2(E)$ be non-zero, and let $K$ be the closure of $\mc B^{(\infty)}(x)$.
As $E$ is essential, it follows that the unit of $\mc A$ acts as the identity on
$E$, and hence also as the identity on $\ell^2(E)$, under $\pi^{(\infty)}$.
Thus $x\in K$.  We shall show that $T(K)\subseteq K$.

Let $V:\ell^2(E)\rightarrow K$ be an $\mc A$-module homomorphism, and let
$\iota:K\rightarrow\ell^2(E)$ be the inclusion map.  By continuity, and the density
of $\mc A$ in $\mc B$, we see that $\iota V \in (\mc B^{(\infty)})'$.
Hence $T \iota V = \iota V T$, from which it follows that $TV(\ell^2(E)) =
VT(\ell^2(E)) \subseteq K$.  Let $W$ be the linear span of the images of all such $V$.
As $\ell^2(E)$ is self-generating, it follows that $W$ is dense in $K$.  However,
$T(W)\subseteq K$, and so by continuity, $T(K)\subseteq K$, as required.

So we have shown that for each $x\in\ell^2(E)$, we have that $T(x)$ is in the
closed linear span of $\mc B^{(\infty)}(x) \subseteq \mc B_w^{(\infty)}(x)$.
By Lemma~\ref{lemma::three}, we conclude that $T \in \mc B_w^{(\infty)}$.
So we have shown that $(\mc B'')^{(\infty)} \subseteq \mc B_w^{(\infty)}$.  By
Lemma~\ref{lemma::one} and Lemma~\ref{lemma::two}, this shows that
$(\mc B^{(\infty)})'' \subseteq (\mc B^{(\infty)})_w$.  However, we always have
that $(\mc B^{(\infty)})_w \subseteq (\mc B^{(\infty)})''$, and so
$(\mc B^{(\infty)})_w = (\mc B^{(\infty)})''$.  Hence also
$(\mc B_w)^{(\infty)} = (\mc B'')^{(\infty)}$, from which it follows immediately
that $\mc B_w = \mc B''$, as required.
\end{proof}

By using the Cohen Factorisation theorem, see \cite[Corollary~2.9.25]{Dales},
a slightly more subtle argument would show that this theorem also holds for
Banach algebras with a bounded approximate identity.

The previous result is only useful if we have a good supply of self-generating
modules.  The following is similar to an idea we used in \cite[Lemma~6.10]{Daws}.

\begin{proposition}\label{prop::one}
Let $\mc A$ be a Banach algebra, and let $E$ be a reflexive Banach space which is
a left $\mc A$-module.  There exists a reflexive left $\mc A$-module $F$ such that:
\begin{enumerate}
\item $E$ is isomorphic to a one-complemented submodule of $F$;
\item each closed, cyclic submodule of $\ell^2(F)$ is isomorphic to
   a one-complemented submodule of $F$;
\end{enumerate}
In particular, $\ell^2(F)$ is self-generating.
\end{proposition}
\begin{proof}
Let $\mc E_0=\{E\}$.  We use transfinite induction to define $\mc E_\alpha$ to be
a set of reflexive left $\mc A$-modules, for each ordinal $\alpha\leq\aleph_1$.
If $\alpha$ is a limit ordinal, we simply define $\mc E_\alpha = \bigcup_{\beta<\alpha}
\mc E_\beta$.

Otherwise, we let $E_\alpha$ to be the $\ell^2$ direct sum of each module in
$\mc E_\alpha$, so that $E_\alpha$ is a reflexive left $\mc A$-module in the obvious
way.  Let $\mc E_{\alpha+1}$ be $\mc E_\alpha$ unioned with the set of all closed cyclic
submodules of $\ell^2(E_\alpha)$.

Let $F$ be the $\ell^2$ direct sum of all the modules in $\mc E_{\aleph_1}$.  As
$\{E\} = \mc E_0 \subseteq \mc E_{\aleph_1}$, condition (1) follows.  Let $K$ be
a closed, cyclic submodule of $\ell^2(F)$, say $K$ is the closure of $\mc A\cdot x$.
Thus
\[ x\in\ell^2(F) \cong \ell^2-\bigoplus_{G\in\mc E_{\aleph_1}} \ell^2(G). \]
Say $x=(x_G)_{G\in\mc E_{\aleph_1}}$ where each $x_G \in \ell^2(G)$.  As
$\|x\|^2 = \sum_G \|x_G\|^2 < \infty$, it follows that $x_G\not=0$ for at most
countably many $G$.  As $\aleph_1$ is uncountable, we must actually have that there
exists $\alpha<\aleph_1$ with $x\in \ell^2-\bigoplus_{G\in\mc E_\alpha} \ell^2(G)
\cong \ell^2(E_\alpha)$.  Then, by construction, $K\in\mc E_{\alpha+1}$, and so $K$
is a one-complemented submodule of $F$.
\end{proof}

Let $\mc A$ be a Banach algebra.  Recall, for example from \cite{Daws}, that
$\wap(\mc A^*)$ is the closed submodule of $\mc A^*$ consisting of those functionals
$\phi\in\mc A^*$ such that
\[ \mc A\rightarrow\mc A^*; \quad a \mapsto a\cdot\phi \]
is weakly-compact.  Young's result, \cite{young}, shows that for each
$\phi\in\wap(\mc AT^*)$, there exists a reflexive Banach space $E$, a contractive
homomorphism $\pi:\mc A\rightarrow\mc B(E)$, and $x\in E,\mu\in E^*$ with
$\|\phi\| = \|x\| \|\mu\|$ and such that
\[ \ip{\phi}{a} = \ip{\mu}{\pi(a)(x)} \qquad (a\in\mc A). \]
Let $\mc A$ be a dual Banach algebra with predual $\mc A_*$.  It is easy to show
(see \cite{Daws} for example) that $\mc A_*\subseteq\wap(\mc A^*)$.  We showed in
\cite[Section~3]{Daws} that Young's result holds for $\phi\in\mc A_*$, with
the additional condition that for any $\lambda\in E^*$ and $y\in E$, the
functional $\pi^*(\lambda\otimes y)$ is in $\mc A_*$, where
\[ \ip{\pi^*(\lambda\otimes y)}{a} = \ip{\lambda}{\pi(a)(y)}
\qquad (a\in\mc A). \]
Note that, a priori, Young's result only shows that $\pi^*(\lambda\otimes y)
\in \wap(\mc A^*)$.

\begin{proposition}\label{prop::two}
With the notation of Proposition~\ref{prop::one}, we have that
$\pi^*(F^*\proten F)$ is a subset of the closed submodule generated by
$\pi^*(E^*\proten E)$.
\end{proposition}
\begin{proof}
The module $F$ is generated from $E$ by two constructions: (i) taking submodules;
and (ii) taking $\ell^2$-direct sums.  For (i), let $K$ be a submodule of
$E$.  The Hahn-Banach theorem shows that $\pi^*(K^*\proten K) \subseteq
\pi^*(E^*\proten E)$.  For (ii), let $(K_i)$ be a family of submodules of $E$ with
$\pi^*(K_i^*\proten K_i) \subseteq \pi^*(E^*\proten E)$ for each $i$, and let
$F = \ell^2-\bigoplus_i K_i$.  Let $\sum_n \mu_n\otimes x_n \in F^*\proten F$,
with, say, $\sum_n \|\mu_n\|^2 = \sum_n \|x_n\|^2 < \infty$.  For each $n$,
we have $\mu_n = (\mu^{(n)}_i)$ with $\|\mu_n\|^2 = \sum_i \| \mu^{(n)}_i \|^2$,
and $x_n = (x^{(n)}_i)$ with $\|x_n\|^2 = \sum_i \| x^{(n)}_i \|^2$.  Then
\[ \sum_n \ip{\mu_n}{a\cdot x_n} = \sum_{n,i} \ip{\mu^{(n)}_i}{a\cdot x^{(n)}_i}
\qquad (a\in\mc A). \]
Hence
\[ \pi^*\Big( \sum_n \mu_n\otimes x_n \Big) = \pi^*\Big( \sum_{n,i} \mu^{(n)}_i
\otimes x^{(n)}_i \Big) \in \pi^*(E^*\otimes E). \]
Again, the Cauchy-Schwarz inequality shows that the sum on the right converges.
\end{proof}

\begin{theorem}\label{thm::two}
Let $\mc A$ be a unital dual Banach algebra.  There exists a reflexive Banach space
$E$ and an isometric, weak$^*$-weak$^*$-continuous homomorphism
$\pi:\mc A\rightarrow\mc B(E)$ such that $\pi(\mc A)'' = \pi(\mc A)$.
\end{theorem}
\begin{proof}
By \cite[Corollary~3.8]{Daws}, we may suppose that $\mc A\subseteq\mc B(E_0)$,
for some reflexive Banach space $E_0$.
By Proposition~\ref{prop::one}, we can find a self-generating, reflexive Banach
space $E$ and a contractive representation $\pi:\mc A\rightarrow \mc B(E)$.
As $E_0\subseteq E$, it follows that $\pi$ is an isometry.  By
Proposition~\ref{prop::two}, $\pi$ is weak$^*$-weak$^*$-continuous.
The result now follows from Theorem~\ref{thm::one}.
\end{proof}

It is well-known that for any Banach algebra $\mc A$, we have that
$\wap(\mc A^*)^*$ is a dual Banach algebra (see, for example,
\cite[Proposition~2.4]{Daws}).  When $\mc A$ has a bounded approximate
identity, a weak$^*$-limit point in $\wap(\mc A^*)^*$ will be a unit for
$\wap(\mc A^*)^*$.

\begin{corollary}
Let $\mc A$ be a Banach algebra with a bounded approximate identity.
There exists a reflexive Banach space
$E$ and a contractive homomorphism $\pi:\mc A\rightarrow\mc B(E)$ such that
$\pi(\mc A)''$ is isometrically, weak$^*$-weak$^*$-continuously isomorphic
to $\wap(\mc A^*)^*$.
\end{corollary}

Finally, we remark that Uygul showed in \cite{uygul} that given a dual, completely
contractive Banach algebra $\mc A$, we can find a reflexive operator space and
a completely isometric, weak$^*$-weak$^*$-continuous homomorphism $\pi:\mc A
\rightarrow\mc B(E)$.  Using this result, we can easily prove a version of
Theorem~\ref{thm::two} for completely contractive Banach algebras.  Indeed, the
only thing to do is to equip $\ell^2$ direct sums with an Operator Space structure
such that the inclusion and projection maps are complete contractions.
This is worked out in detail in \cite{xu} (see also \cite{uygul}).

Finally, we remark that the space constructed in Theorem~\ref{thm::two} is
very abstract.  For a group measure space convolution algebra $M(G)$, Young showed
in \cite{young} that $M(G)$ can be weak$^*$-represented on a direct sum of
$L^p(G)$ spaces; the analogous result for the Fourier algebra was shown by the
author in \cite{daws1}.  For such concrete Banach algebras $\mc A$, it would be
interesting to know if ``nice'' reflexive Banach spaces $E$ could be found with
$\pi:\mc A\rightarrow\mc B(E)$ such that $\pi(\mc A)''=\pi(\mc A)$.

\vspace{5ex}

\noindent\emph{Author's Address:}
\parbox[t]{3in}{School of Mathematics,\\
University of Leeds\\
Leeds\\
LS2 9JT.}

\bigskip\noindent\emph{Email:} \texttt{matt.daws@cantab.net}

\end{document}